\documentclass[11pt]{article}
\usepackage[T1]{fontenc}
\usepackage{microtype}
\usepackage[a4paper,margin=2.5cm]{geometry}
\usepackage{amsmath, amssymb, amsthm}

\usepackage{hyperref, cleveref}

\newtheorem{thr}{Theorem}[section]
\newtheorem{lem}[thr]{Lemma}
\newtheorem{prop}[thr]{Proposition}

\newcommand{\abs}[1]{\lvert #1\rvert}
\newcommand{\floor}[1]{\lfloor #1\rfloor}

\title{The majority game chromatic number for dense random graphs}

\date{}
\author{Stijn Cambie\thanks{Department of Computer Science, KU Leuven Campus Kulak-Kortrijk, 8500 Kortrijk, Belgium. E-mail: \protect\href{mailto:stijn.cambie@hotmail.com}{\protect\nolinkurl{stijn.cambie@hotmail.com}}. Supported by FWO grants 1225224N and 1215727N.} 
\and
    Annika Heckel\thanks{Matematiska institutionen, Uppsala universitet, Box 480, 751 06 Uppsala, Sweden. Email: \href{mailto:annika.heckel@math.uu.se}{\nolinkurl{annika.heckel@math.uu.se}}. Funded by the Swedish Research Council, Starting Grant 2022-02829.}
\and
    Brahadeesh Sankarnarayanan\thanks{Department of Mathematics, Indian Institute of Technology Jodhpur, Jodhpur 342030, Rajasthan, India. Email: \href{mailto:brahadeesh@iitj.ac.in}{\nolinkurl{brahadeesh@iitj.ac.in}}. Supported by the Research Initiation Grant I/RIG/BDS/20250230, IIT Jodhpur.}
\and
Gabri\"elle Zwaneveld\thanks{Korteweg-de Vries Instituut, Universiteit van Amsterdam, Sciencepark 107, 1098 XG Amsterdam. Email: \href{mailto:g.c.zwaneveld@uva.nl}{\nolinkurl{g.c.zwaneveld@uva.nl}}}
}

\begin{document}
\maketitle


\begin{abstract}
    We show that for every fixed $2/3<p<1$ the majority game chromatic number of the Erd\H{o}s--R\'enyi random graph $G_{n,p}$ is $3$ with high probability. 
\end{abstract}

\section{Introduction}

Let \(G = (V,E)\) be a finite, simple graph.
A (not necessarily proper) vertex colouring of \(G\) is said to satisfy the \emph{majority condition} at \(v \in V\) if at most half the neighbours of \(v\) have the same colour as $v$.
The \emph{majority colouring game}, introduced by Bosek--Grytczuk--Jak\'obczak~\cite{BGJ19} in 2019, is a two-player Maker-Breaker-type game in which the players alternately colour the vertices of \(G\) from a fixed palette of colours such that the majority condition is maintained at every vertex at all times.
The objective of the first player, Alice, is to get a colouring of all the vertices, which is called a \emph{majority colouring} of the graph \(G\).
The objective of the second player, Bob, is to prevent Alice from succeeding; that is, Bob tries to create a partial colouring where there is an uncoloured vertex for which every choice of colour violates the majority condition somewhere in the graph.
We assume that Alice moves first.
The minimum number of colours required for Alice to have a winning strategy for the majority colouring game on \(G\) is called the \emph{majority game chromatic number} of \(G\) and is denoted \(\mu_g(G)\).

In~\cite{BGJ19}, the authors proved that \(\mu_g(G)\) is unbounded on the class of bipartite graphs and on the class of \(2\)-degenerate graphs. They asked whether $\mu_g(G_{n,p})$ is bounded for random graphs (which we answer in Proposition~\ref{Prop:ub}), and conjectured that it is usually two~\cite[Conj.~2]{BGJ19}.
In this short note (in Theorem~\ref{thm:main}), we prove that $\mu_g(G_{n,p})$ is usually three for a dense random graph \(G_{n,p}\). We also make a number of quick observations in Proposition~\ref{prop:observations} which bound or determine $\mu_g(G_{n,p})$ for other ranges of $p=p(n)$.

\subsection*{Notations and preliminaries}

The (Erd\H{o}s--R\'enyi) random graph $G_{n,p}$ is a probability distribution over all graphs obtained on a vertex set of order $n$, with every possible edge independently sampled with probability $p.$
The (open) neighbourhood of a vertex $v$ is denoted with $N(v)$.
The degree of $v$ is $d(v)=|N(v)|$.
The minimum degree of a graph $G$ is denoted by $\delta(G)$.
The symmetric difference of two sets is denoted as $A \triangle B= (A \setminus B) \cup (B \setminus A).$
A sequence $(E_n)_{n \in \mathbb{N}}$ of events holds \emph{with high probability (whp)} if $\Pr(E_n) \rightarrow 1$ as $n \rightarrow \infty$.
Throughout, we use the known consequences of Chernoff's bounds (see \cite[Cor.~2.3]{janson:randomgraphs}) on the concentration of various graph parameters of $G_{n,p}.$
\begin{lem}
\label{chernoffcorollary}
Let $X$ be distributed binomially with parameters $n $ and $p $. If $0 <\epsilon \leq 3/2$, then
\[
\Pr\left( \left| X - np \right| \geq \epsilon np \right) \leq 2 e^{-\epsilon^2 np /3}.
\]
\end{lem}

\section{Bounds for random graphs}
First, we give an elementary proof that the majority game chromatic number is whp bounded for $G_{n,p}$, where $p>0$ is a constant.

\begin{prop}\label{Prop:ub}
    Let $p>0$ be a constant. Then, whp, $\mu_g(G_{n,p})\le \floor{2/p}+1$.
\end{prop}

\begin{proof}
    Let $r=\floor{2/p}+1.$ 
    By Chernoff's bound~\ref{chernoffcorollary}, the probability that the minimum degree of $G \sim G_{n,p}$ is less than $\frac{2n}{r}$ is exponentially small.
    So, assume $\delta(G)\geq\frac{2n}{r}$, and consider a partial colouring of \(G\) in the majority colouring game on \(G\) with \(r\) colours, where Alice plays any legal move on her turn regardless of Bob's strategy.
    A least common colour appears at most 
    $\frac{n}{r}\leq\frac{\delta(G)}{2}$ many times, so such a colour can be played on any uncoloured vertex.
    Thus both players always have a legal move available on their turn, and hence Alice wins.
\end{proof}

As a corollary, this implies that $\mu_g(G_{n,p}) \le 3$ whp when $p>2/3$ is a constant.
Next, we disprove~\cite[Conj.~2]{BGJ19} for $G_{n,p}$ with $2/3<p<1$ a constant, by proving that whp $\mu_g(G_{n,p})= 3$. We conjecture that $\mu_g(G_{n,p}) = 3$ whp for every fixed $0<p<1$.

\begin{thr}\label{thm:main}
    If $2/3<p<1$ is a constant, then for $G \sim G_{n,p}$, whp $\mu_g(G) \ge 3$.
\end{thr}

\begin{proof}
    We prove that for such $G \sim G_{n,p}$, whp Bob has a winning strategy for the majority colouring game played with two colours (which we call red and blue).
    Without loss of generality, assume Alice first colours a vertex $v$ in red.
    Next, Bob colours a neighbour $w$ of $v$ in red, and then Alice colours a vertex $u$.

    After this, Bob colours a vertex $z$ in blue, where $z$ satisfies the following conditions
    \begin{itemize}
        \item it is a non-neighbour of both $v$ and $w$,
        \item $u$ is a neighbour of $z$ if and only if $u$ is a neighbour of $v$,
        \item $\abs{N(v)}+\abs{N(z)}$ is odd, and
        \item $3\abs{N(v) \cap N(z)} > \abs{N(v)}+ \abs{N(z)}$ (equivalently, $\abs{N(v) \cap N(z)}> \abs{N(v) \triangle N(z)}$).
    \end{itemize}

    Writing \(I = N(v) \cap N(z)\), \(A = N(v) \setminus N(z)\) and \(B = N(z) \setminus N(v)\), the last condition is that \(|I|>|A|+|B|\), so that \(v\) and \(z\) have many common neighbours.
Such a vertex \(z\) exists with high probability, since the number of vertices satisfying the first $3$ conditions is at least $0.5(1-o(1))n (1-p)^3$ and almost all of them satisfy the last property, since $3p^2> 2p.$ 

In a majority colouring, \(v\) has at most \(\floor{d(v)/2}\) red neighbours and \(z\) has at most \(\floor{d(z)/2}\) blue neighbours.
Every vertex of \(I\) uses up one unit of this combined budget, whatever its colour.
Since \(d(v) + d(z)\) is odd, the budget left over for \(A \cup B\) is \((|A| + |B| - 1)/2\).
So if more than half of the vertices of \(A \cup B = N(v) \triangle N(z)\) receive the same colour as their neighbour in \(\{v, z\}\), the position has no legal completion, and Bob wins.

In the next \((|A| + |B| - 1)/2\) moves of Bob, he colours a vertex in $N(z) \triangle N(v)$ in the same colour as the neighbour in $\{v,z\}$, in a balanced way, i.e., such that in $A$ and $B$ there are respectively at least $|A|/2-2$ and $|B|/2-2$ many vertices coloured by him.
To ensure this is possible, as long as fewer than the desired number of vertices of $A$ (resp., $B$) have been coloured by Bob, when Alice colours a vertex in $A$ (resp., $B$), Bob chooses to do so in that part as well.
Otherwise, one of the two parts $A$ and $B$ for which the desired number of coloured vertices is not obtained can be freely chosen. 
    Since whp $|A|, |B| \sim np(1-p) $, no more than $(1+o(1)) \frac 32 np(1-p)< \frac 12 \delta(G)$ vertices were coloured with the same colour and thus no forbidden move was made at this point.    
    
    In any case, after these turns, more than half of the vertices of $N(z) \triangle N(v)$ have the same colour as their neighbour in $\{v,z\}$.
    So, if the game ended with a majority colouring of \(G\), then either $v$ or $z$ would have more than half of its neighbours in its own colour, a contradiction.
    Thus Bob wins and $\mu_g(G)>2$.
\end{proof}

Lastly, we prove some results about \(\mu_g(G_{n,p})\) when \(p\) is not constant.
\begin{prop}  \label{prop:observations}
The following statements all hold with high probability for $G \sim G_{n,p}$.
\begin{enumerate}
\item If \(n^2 p \to 0\), then \(\mu_g(G) = 1\).
\item If \(n^{-2} \ll p \ll n^{-7/6}\), then \(\mu_g(G) = 2\). 
\item If \(n^{-7/6} \ll p \le C/n\) for a fixed $C > 0$, then \(\mu_g(G) \ge 3\).
\item If \(p = 1 - o(n^{-3/2})\), then \(\mu_g(G) = 2\).
\item If \(p = 1 - 1/n\) and \(2 \nmid n\), then \(\mu_g(G) = 3\).
\end{enumerate}
\end{prop}

\begin{proof}
\begin{enumerate}
\item If \(n^2 p \to 0\), then whp \(G\) is empty.
\item When \(n^{-2} \ll p \ll n^{-7/6}\), whp each component of \(G\) is a tree on at most six vertices, while whp $G$ is not empty.
An exhaustive check shows that each of them is an Alice-win with \(2\) colours, regardless of which player goes first and even if Bob skips some of his turns.
\item Bob wins on a palette of \(2\) colours as soon as \(G\) has two components isomorphic to the \(1\)-subdivision of \(K_{1,3}\) (see~\cite{CNSS2026-CCHLSYZ}), so whp \(\mu_g(G) > 2\) when \(n^{-7/6} \ll p \le C/n\).
\item When $p=1-o(1/ n^{3/2})$, whp the complement of \(G\) is a matching.
    If $n$ is odd, Alice can ensure that non-neighbours have the same colour.
    If $n$ is even, then \(\delta(G) \geq n - 2\), so a least common colour is used on at most \(\floor{(n-1)/2} \leq \delta(G)/2\) vertices. So, as in the proof of Proposition~\ref{Prop:ub}, a legal move always exists.

\item For $p=1-1/n$, $n$ odd, whp there are at least $3$ components of the complement that are cherries, i.e., paths on three vertices. If Alice and Bob play the game with two colours, Bob ensures in the early phase of the game that there are two cherries that are not monochromatic and whose centres have different colours.
    Suppose \(G\) ends up with a majority colouring.
    Then a non-monochromatic cherry whose centre is not in the larger colour class has a leaf in the larger colour class that violates the majority condition, a contradiction.
    So, \(\mu_g(G) \geq 3\).
    A similar argument as in the proof of Proposition~\ref{Prop:ub} shows that \(\mu_g(G) \leq 3\).
\end{enumerate}
\end{proof}

\subsection*{Acknowledgements}

This problem was proposed during the MATRIX--MFO Tandem Workshop: Combinatorial Interchange, in September 2026.

\subsection*{AI declaration}
The main result and its underlying proof strategy were obtained by the authors. Generative AI was used to check draft arguments, and to explore possible extensions. Those extensions are not claimed as results in this paper. In particular, AI claimed to have proved that $\mu_g(G_{n,p})=3$ whp for every $0<p<1.$ We only conjecture the correctness of the latter extension here.

\end{document}